\documentclass[twoside,a4paper,preprint]{imsart}
\usepackage{subcaption}
\usepackage{pdflscape}
\usepackage{float}
\usepackage[breaklinks=true]{hyperref}
\usepackage{breakcites}
\usepackage[utf8]{inputenc}
\usepackage{bbold}

\usepackage{booktabs}
\usepackage{multirow,color}
\usepackage{rotating}
\usepackage{diagbox}

\usepackage[shortlabels]{enumitem}
\usepackage{bm}

\usepackage{xcolor-patch}

\usepackage{amsmath}
\usepackage{amssymb}
\usepackage{amsthm}

\usepackage{boxedminipage}

\usepackage{listings}
\usepackage{minitoc}
\usepackage{natbib}

\usepackage{ifpdf}
\usepackage{lipsum}

\DeclareMathOperator{\logit}{logit}

\newcommand{\1}{OR$_{X_Y}$}
\newcommand{\2}{AIPW$_{X_Y}$}
\newcommand{\3}{OR$_{DS}$}
\newcommand{\4}{AIPW$_{X_Y,X_T}$}
\newcommand{\5}{AIPW$_{DS}$}
\newcommand{\6}{TMLE$ _{X_Y, X_T}$}

\newcommand{\ci}{\perp\!\!\!\perp}
\theoremstyle{definition}
\newtheorem{theorem}{Theorem}

\newtheorem{corollary}{Corollary}

\newtheorem{assumption}{Assumption}

\RequirePackage[OT1]{fontenc}
\usepackage{amsthm,amsmath,natbib}
\RequirePackage[colorlinks,citecolor=blue,urlcolor=blue]{hyperref}

\startlocaldefs
\numberwithin{equation}{section}
\theoremstyle{plain}

\endlocaldefs

\begin{document}

\begin{frontmatter}
\title{The costs and benefits of uniformly valid causal inference with high-dimensional nuisance parameters}
\runtitle{UNIFORMLY VALID CAUSAL INFERENCE}

\author{\fnms{Niloofar} \snm{Moosavi}},
\author{\fnms{Jenny} \snm{H\"aggstr\"om}}
\and
\author{\fnms{Xavier} \snm{de Luna}}
\runauthor{N MOOSAVI, J H\"AGGSTR\"OM AND X DE LUNA}

\affiliation{Umeå School of Business, Economics and Statistics at Umeå University}
\address{Department of Statistics, USBE, Ume\aa{} University, 901 87, Ume\aa, Sweden}
\address{niloofar.moosavi@umu.se; jenny.haggstrom@umu.se; xavier.deluna@umu.se}

\begin{abstract}
Important advances have recently been achieved in developing procedures yielding uniformly valid inference for a low dimensional causal parameter when high-dimensional nuisance models must be estimated. In this paper, we review the literature on uniformly valid causal inference and discuss the costs and benefits of using uniformly valid inference procedures. Naive estimation strategies based on regularisation, machine learning, or a preliminary model selection stage for the nuisance models have finite sample distributions which are badly approximated by their asymptotic distributions. To solve this serious problem, estimators which converge uniformly in distribution over a class of data generating mechanisms have been proposed in the literature. In order to obtain uniformly valid results in high-dimensional situations, sparsity conditions for the nuisance models need typically to be made, although a double robustness property holds, whereby if one of the nuisance model is more sparse, the other nuisance model is allowed to be less sparse. While uniformly valid inference is a highly desirable property, uniformly valid procedures pay a high price in terms of inflated variability. Our discussion of this dilemma is illustrated by the study of a double-selection outcome regression estimator, which we show is uniformly asymptotically unbiased, but is less variable than uniformly valid estimators in the numerical experiments conducted. 
\end{abstract}


\end{frontmatter}

\section{Introduction}
High-dimensional situations, where the number of covariates is larger than the number of observations are common in causal inference applications. Using regularization type estimators such as lasso \citep{tibshirani1996regression} or other post-model selection estimators are popular strategies in such cases. Important advances have been achieved in developing procedures yielding uniformly valid inference (defined below) for a low dimensional causal parameter when high-dimensional nuisance models must be fitted \citep[e.g.,][]{van2006targeted,belloni2014inference,van2014targeted,farrell2015robust,chernozhukov2018double}.  In this paper, we review the literature on uniformly valid causal inference, and discuss the costs and benefits of using uniformly valid inference procedures. This discussion is important since naive and invalid post-model selection inference is to this day still common in statistical practice.

\citet{leeb2005model} demonstrated how a data-driven model selection step
can affect the distribution of the estimate of
a parameter of interest. 
 Loosely, they show that the scaled  ($\sqrt{n}$) bias  of a naive two step estimator, which does not take into account the selection step, goes to infinity or stay bounded for a sequence of worst case scenario data generating processes (DGPs), when relying on consistent or conservative model selection, respectively. We say that such a naive estimator is not uniformly unbiased. An estimator with associated uniformly valid inference, on the other hand, is such that its distribution $F_n$ converges uniformly over a set of DGPs $\mathcal{P}$, i.e., for any $u \in \mathbb{R}$
 $$
 \lim_{n \rightarrow \infty} \sup\limits_{P \in \mathcal{P}}|F_n(u) - F_P(u)| =0,
 $$
where $F_P$ is the cumulative distribution derived from $P$. Hereafter, we use the terms uniformly valid inference and valid inference interchangeably. 

Valid inference for every single parameter in a linear regression model needs careful consideration in high-dimensional settings. Some have considered debiasing lasso for a valid inference that target the true data generating process parameters \citep{javanmard2014confidence,van2014asymptotically,zhang2014confidence} and others have considered valid inference conditional on the model that has been selected \citep{berk2013valid,lee2016exact}. In a causal inference context, there is typically a low dimensional parameter of interest, e.g., the average causal effect of a treatment, and high-dimensional nuisance parameters. \citet{belloni2014inference} proposed an estimator with valid inference for a causal parameter in a linear model explaining outcome with a treatment variable and a set of covariates, which can be of high dimension. To achieve uniformly valid inference, they proposed to include the union of two sets of covariates in the model: one obtained by selecting covariates relevant when regressing the outcome on the covariates, and the second by selecting covariates relevant when regressing the treatment on the covariates. Their model implicitly implies a homogeneous causal effect. This can be relaxed to allow for individual heterogeneous effects, using the potential outcome framework \citep{neyman1923application,rubin1974estimating}, and nuisance models for both the potential outcomes and for the treatment assignment given the covariates (propensity score). For this general case, \citet{van2006targeted} and \citet{van2014targeted} obtained valid inference for a causal parameter using targeted maximum likelihood estimation, where nuisance models are estimated nonparametrically. \citet{farrell2015robust} considered the augmented inverse probability of treatment weighting estimator \citep{robins1994estimation}, and showed that uniformly valid inference is achieved when using post-lasso estimation for the nuisance models. Similar results were derived in \cite{chernozhukov2018double} using a double machine learning approach.

In both \citet{belloni2014inference} and \citet{farrell2015robust}, approximate sparsity is assumed for the nuisance models. However, in \citet{farrell2015robust} the outcome model can be less sparse if the propensity score is more sparse and vice versa \citep[called nonparametric double robustness property,][]{kennedy2016semiparametric}. Yet consistency in the nonparametric estimation of all the nuisance models is required, a condition relaxed in \citet{van2014targeted}, and in more recent work \citep{avagyan2017honest,tan2020model}, where one of the nuisance models may be inconsistently estimated.

Procedures yielding uniformly valid inference for a causal parameter, in the general context of heterogeneous treatment effects, allow for the selection of instruments (loosely, variables related to the treatment but not the outcome) in the fit of the propensity score model. This is known to result in possibly large inflation of the variance of the estimators \citep[e.g.,][]{hahn2004functional,de2011covariate,Schnitzer2016,RotnitzkySmucler:2020}. Thus, while uniformly valid inference is a highly desirable property, uniformly valid procedures pay a high price in terms of inflated variability. We discuss and illustrate this dilemma by studying a compromise solution, an outcome regression estimator \citep[e.g.,][]{Tan:2007}, which we allow to select instruments, but which does not use the fitted propensity score, in contrast with uniformly valid estimators proposed in the literature. The resulting post-model selection estimator is shown to be uniformly asymptotically unbiased under a commonly used product rate condition \citep{farrell2015robust}, even though the propensity score is not used in the estimator except for the covariate selection step. 

This paper is organised as follows. Section 2 presents a review of the literature on uniformly valid causal inference. Section 3 gives a theoretical discussion of the costs and benefits of uniformly valid inference, by studying a double-selection outcome regression estimator. Section 4 illustrates this discussion with a Monte Carlo study of finite sample properties of a collection of estimators. Section 5 concludes the paper. All proofs are delayed to an appendix.

\section{Uniformly valid causal inference: a review}\label{uniformlyvalidinference}
This is an extremely active research area. The focus is here on uniformly valid inference on a low dimensional causal parameter after regularization/model selection of high-dimensional nuisance models. We start by introducing some general concepts, and then review first cornerstone work on the homogeneous and then general heterogeneous case. A review of important advances in recent years concludes this section.

The parameter of interest is a causal effect of a binary treatment variable $T$ on an outcome $Y$ as defined below in different contexts. We use the notation $X$ to denote a one dimensional pretreatment covariate and $\bm{X}$ to denote a set of pretreatment covariates which has dimension $p$, allowed to grow with $n$. Note that the set may contain not only the covariates but also transformations of them. We consider a set of identically and independently distributed (i.i.d.) observations, $\{(\bm{x}_i, y_i, t_i)\}^{n}_{i=1}$, drawn from a distribution $P_n$. To study uniformly valid post-model selection inference it is essential that the probability law $P_n$ is allowed to vary with the sample size $n$. We use further the notations $E_n[w_i ] = \frac{1}{n}\Sigma_{i=1}^n w_i$ and $a \vee b = \max\{a, b \}$. Moreover, $n_t$ denotes the number of individuals under treatment.

\subsection{Hodges estimator and superefficiency}
Superefficient estimators of a parameter of a model $\mathcal{M}$ are variants of the well known Hodges estimator \citep{vandervaart:97} when a model restriction holds; meaning that a model $\mathcal{M}_0 \subset \mathcal{M}$ contains the true data generating process. The asymptotic variance of a superefficient estimator is smaller than the efficiency bound for the class of regular asymptotic linear (RAL) estimators of the parameter under model $\mathcal{M}$, when there is a submodel $\mathcal{M}_0$ under which the same bound is smaller. Superefficiency has a cost in the sense that the asymptotic distribution of a superefficient estimator is valid only pointwise at $\mathcal{M}_0$ instead of uniformly over a larger family of models (uniformly valid inference). 

This was highlighted by \citet{leeb2005model} in a parametric setting (see below for a detailed exposition), where a consistent model selection step that selects out a "redundant" variable (not part of $\mathcal{M}_0$) before a maximum likelihood fit results in a superefficent estimator. Such an estimator has an oracle property in the sense that it asymptotically (only pointwise at $\mathcal{M}_0$ instead of uniformly over a larger family of models) performs as well as a fictitious orcale estimator which can be constructed by knowing $\mathcal{M}_0$ \citep{fan2001variable}.  

\subsection{Homogenous causal effect}\label{homog.sec}

As a primer, consider the parametric regression model $y_i = \alpha t_i + \beta x_i +\epsilon_i$, with $\epsilon_i \sim  \mathcal{N}(0, \sigma^2),\, \sigma^2 > 0$.  
According to \citet{leeb2005model}, if we are interested in $\alpha$,
 the post-selection estimator which includes a preliminary consistent model selection step on $X$ (i.e., a test whether $\beta =0$) is more efficient than the simple OLS estimator without this step if corr$(X,T)\neq 0$. However, if $\beta\neq 0$ and corr$(X,T)\neq 0$, the finite sample distribution of the post-selection estimate is a mixture of two normal distributions, i.e., not well approximated by the normal asymptotic distribution. This is because in the selection step the nonzero coefficient can be detected for some samples and not detected for others. If the resulting omitted variable bias is considerable, the empirical coverage of a naive confidence interval can be far from the nominal coverage. \citet{leeb2005model} have shown that the minimal coverage of the naive confidence interval with respect to all possible $\beta$ values goes to zero as $n$ grows for consistent model selection steps, while the empirical coverage for any fixed $\beta$ value goes to the nominal one. Their result highlights the importance of uniformly valid inference compared to pointwise asymptotic results.
 
 For $\alpha$ to have a causal interpretation, the linear regression needs to include all confounders as formally defined in next section, and it must be correctly specified as a model for $E(Y\mid T,\bm{X})$
(in particular implying a homogeneous/constant causal effect). Concerning the former condition, the number of available covariates, hence potential confounders, may be very large when using large observational databases. The latter condition, implicitly requires that series expansions need to be used to approximate $E(Y\mid T,\bm{X})$ increasingly well with increasing sample sizes. Therefore, these two conditions often yield a high-dimensional setting in practice, i.e. where the number of covariates is at least as large as the sample size.

 \citet{belloni2014inference} proposed a strategy for reaching valid inference in such high-dimensional settings \citep{LEEB2008}. Their suggestion is a two-step lasso-based method, where the union of covariate sets selected by two distinct lasso regressions of $Y$ and $T$ on the covariates, respectively (often called double selection), are utilized in a second step in the main linear model including $T$ as a regressor. Instead of exact sparsity conditions typically used in high-dimensional settings, they consider the following approximate sparsity conditions.
Let 
\begin{equation}\label{bellonimodely}
E(Y|T,\bm{X}) = \alpha T + \beta_Y' \bm{X} + R_Y,
\end{equation}
\begin{equation}\label{bellonimodelt}
E(T|\bm{X}) =\beta_T' \bm{X} + R_T,
\end{equation}
where $R_f$ is the specification error of using a sparse $\beta_f$ with only $s_f$ nonzero elements, respectively for $f=Y, T$. The regularity conditions to obtain uniform valid inference include
$$E(E_n[R_{f,i}^2])^{1/2} = O\big(\sqrt{s_f/n}\big),$$
$$\log^3 p/n =o(1),$$
$$s_f^2\log^2(p \vee n)/n = o(1).$$ 
In other words, models (\ref{bellonimodely}) and (\ref{bellonimodelt}) are assumed well approximated  by a sparse linear combination of the covariate vector $\bm{X}$. 
 Under these conditions (and other regularity conditions), \citet[][Corollary 2]{belloni2014inference} showed that their estimator of $\alpha$ is asymptotically normal uniformly over $P_n$, thereby uniformly valid inference can be made.

\subsection{Heterogeneous causal effect}\label{hetero.sec}
The effect of a binary treatment is now allowed to be heterogeneous using the Neyman-Rubin potential outcome framework \citep{neyman1923application,rubin1974estimating}. For any unit in the study, denote $Y(1)$ its potential outcome under treatment ($T=1$),  and $Y(0)$ its potential outcome without treatment (or alternative treatment). We assume that $Y = TY(1) + (1-T)Y(0)$ is the observed outcome, and no interference between units is allowed \citep[stable unit treatment value assumption;][]{DR:90}.
Each unit may have a different causal effect Y(1) - Y(0), and the average causal effect $\tau = E(Y(1) - Y(0))$ is the parameter of interest in the sequel. This parameter is identified given the following assumptions.
\begin{assumption}[No unobserved confounding]\label{NUC}
$Y(1), Y(0)\ci{T|\bm{X}}. $
\end{assumption}
\begin{assumption}[Overlap]\label{overlap}
$\mathbb{P}(T=t|\bm{X}) \geq p_{min} > 0 ,  \qquad t=0, 1.$
\end{assumption}
\noindent Thus, all confounding covariates are included in $\bm{X}$ and all units in the study have non-zero probability to be included in both treatment groups. These assumptions are made throughout the article. Note, however, that a sensitivity analysis \citep[see][and references therein]{diazetal18} should accompany inference based on Assumption 1, even when many covariates are available, since this assumption is not testable without further information \citep[e.g.,][]{LunaJ:2014}. Assumption 2 has also important implications as was recently demonstrated in \citet{DAMOUR:21} for high-dimensional situations, where overlap is linked to the sparsity conditions discussed herein.
 
 Let $E(Y|T=1, \bm{X}) = m_1(\bm{X})$ and $E(Y|T=0, \bm{X}) = m_0(\bm{X})$ denote the outcome models, $E(T|\bm{X}) = \mathbb{P}(T = 1|\bm{X})= e(\bm{X})$ denotes the propensity score model. One of the earliest proposals that addressed inference on $\tau$ when estimating nuisance models nonparametrically is the targeted maximum likelihood estimator \citep[TMLE;][]{van2006targeted, vdLR:11}. Denote fits of the nuisance models $\mathbb{P}(T=1|\bm{X})$, $E(Y(1)|\bm{X})$ and $E(Y(0)|\bm{X})$ by $\hat{e}(\bm{x}_i)$, $\hat{m}^0_1(\bm{x}_i)$ and $\hat{m}^0_0(\bm{x}_i)$, respectively. Then, a fluctuated version of the predicted outcome values is used in the following manner:
\begin{equation*} 
\label{tmle}
     \hat{\tau}_{TMLE} = E_n [ \hat{m}^1_1(\bm{x}_i) -
   \hat{m}^1_0(\bm{x}_i)],
\end{equation*}
where the fluctuations are found by
\begin{equation*}
    \logit \hat{m}^1_t(\bm{x}_i) =  \logit \hat{m}^0_t(\bm{x}_i) + \varepsilon_n  h_t(\bm{x}_i) , \qquad t \in \{0,1\},
\end{equation*}
where
\begin{equation*}
    h_t(\bm{x}_i) = \dfrac{\mathbb{1}\{t=1\}}{\hat{e}(\bm{x}_i)} - \dfrac{\mathbb{1}\{t=0\}}{1 - \hat{e}(\bm{x}_i)},
\end{equation*}
and $\varepsilon_n$ is found by running logistic regression of outcome $Y$ on $h_T(\bm{X})$ using $\logit \hat{m}^0_T(\bm{X})$ as intercept.
TMLE is consistent if either $e(\cdot)$ or $m_0(\cdot)$ and $m_1(\cdot)$ are consistently estimated. Moreover, it is RAL and semiparametrically efficient if all models are consistently estimated and a product of rate of convergence similar to (\ref{productrate.eq}) hold.
TMLE can also be constructed by iteratively fluctuating the propensity score and the outcome models, thereby yielding a RAL estimator when at most one of the models is consistently estimated, i.e. this TMLE is then not only consistent but also asymptotically normal
\citep[so called double robust statistical inference,][]{van2014targeted}.

In another major contribution,
\citet{farrell2015robust} showed how to obtain uniformly valid inference for the popular doubly robust augmented inverse probability of treatment weighting \citep[AIPW,][]{robins1994estimation, scharfstein1999rejoinder} estimator:
\begin{equation} 
\label{dr}
     \hat{\tau}_{DR} = E_n \left[     \dfrac{t_iy_i - (t_i - \hat{e}(\bm{x}_i))\hat{m}_1(\bm{x}_i)}{\hat{e}(\bm{x}_i)} -
    \dfrac{(1-t_i)y_i + (t_i - \hat{e}(\bm{x}_i))\hat{m}_0(\bm{x}_i)}{1-\hat{e}(\bm{x}_i)}
    \right].
\end{equation}

\noindent The fitted values $\hat{m}_0(\bm{x}_i)$, $\hat{m}_1(\bm{x}_i)$ and $\hat{e}(\bm{x}_i)$ are obtained using, e.g., post-lasso estimators. \citet{farrell2015robust} proposed the use of group-lasso to benefit from grouped sparsity patterns among potential outcomes and different treatment levels. Similar to \citet{belloni2014inference}, \citet{farrell2015robust} assumed approximate sparsity, but for logistic propensity score and linear potential outcome models, i.e.,
\begin{equation}\label{farrellmodely}
m_t(\bm{X},\eta_t)= \eta'_{t}\bm{X}+ R_Y^t , \qquad t \in \{0,1\}, 
\end{equation}
\begin{equation*}\label{farrellmodelt}
e(\bm{X},\gamma)= \text{expit}(\gamma'\bm{X} + R_T),
\end{equation*}
where $R_Y^t$ and $R_T$ are approximation errors of estimating the true models with sparse parameters $\eta_t$ and $\gamma$ that have $s_Y^t$ and $s_T$ nonzero elements, respectively. Farrell's regularity conditions on the specification errors are slightly different from those in \citet{belloni2014inference}. In particular, he assumes
\begin{equation}\label{errorY}
    E_n[(r_{Y,i}^t)^2]^{1/2} \vee E_{n_t}[(r_{Y,i}^t)^2]^{1/2} \leq \mathcal{R}_Y^t= O(\sqrt{s_Y^t/n}), \qquad t \in \{0,1\},
\end{equation}
and
\begin{equation}\label{errorT}
    E_n[(\text{expit}(\gamma'\bm{x}_i)-\text{expit}(\gamma'\bm{x}_i + r_{T,i}))^2]^{1/2} \leq \mathcal{R}_T = O(\sqrt{s_T/n}).
\end{equation}
More importantly, the sparsity assumption required for each nuisance model separately is weaker compared to the one assumed by \citet{belloni2014inference}, i.e., 
\begin{equation}\label{consistencysparsity}
    s_f \log(p \vee n )^{3/2+\delta}
=o(n), \qquad s_f \in \{s_T, s_Y^0, s_Y^1\},
\end{equation}
for some $\delta>0$, since here we have a multiplicative rate condition 
\begin{equation}\label{ratesparsity}s_Y^t s_T \log(p \vee n )^{3+2\delta}=o(n), \qquad t \in \{0, 1 \}.\end{equation}
Thus, if the potential outcome models are more sparse, the propensity score model is allowed to be more dense and the other way around. These sparsity assumptions and other regularity conditions result in the following rates of convergence for 
the post-lasso estimators of the nuisance models \citep[][Section 6]{farrell2015robust}
$$E_n[(\hat{e}(\bm{x}_i) - e(\bm{x}_i))^2] = o_{P_n}(1), \quad E_n[(\hat{m}_t(\bm{x}_i) - m_t(\bm{x}_i))^2] = o_{P_n}(1),$$
\begin{align}\label{productrate.eq}
    E_n[(\hat{e}(\bm{x}_i) - e(\bm{x}_i))^2]^{1/2} E_n[(\hat{m}_t(\bm{x}_i) - m_t(\bm{x}_i))^2]^{1/2}=o_{P_n}(n^{-1/2}),
\end{align}
for $t \in \{0,1\}$. Under these consistency and product rate conditions (and other regularity assumptions) the estimator (\ref{dr}) is asymptotically normal uniformly over $P_n$ \citep[][Corollary 3]{farrell2015robust}. This result is not restricted to the post-lasso estimator, but ensures $\sqrt{n}$-consistency of the AIPW estimator of the low dimensional parameter of interest $\tau$ for any estimator of nuisance models which fulfills the assumptions. High-dimensional parametric or nonparametric nuisance models fit into this framework even though the estimation cannot be done at the $\sqrt{n}$-rate.

The presentation above has focused on robustness to the danger of including too few covariates (regularization bias), in settings where we believe in sparsity assumptions. Another possible source of bias arises from the danger of overfitting. This is a problem when the nuisance functions are too complex (e.g., cannot be assumed to belong to a Donsker class; e.g., \citeauthor{diaz:2019}, \citeyear{diaz:2019};  \citeauthor{kennedy2016semiparametric}, \citeyear{kennedy2016semiparametric}). A general solution to avoid overfitting error, and thereby obtain valid inference, is to use sample-splitting; see, \cite{Bickeln:82,zheng2011cross,cattaneoetal:2018,cattaneoetal:2019}.

\subsection{Advances in recent years}
Recent years have witnessed a great deal of novel results in the field of uniformly valid causal inference. \cite{chernozhukov2018double} readdressed 
valid inference for the average causal effect under the framework of double/debiased machine learners in light of the fact that the parameter $\tau$ satisfies a Neyman orthogonal moment condition \citep{Neyman:59,Neyman:79}; a moment condition that is not sensitive to local errors in nuisance models and can be derived using the first order influence function of the parameter \citep{bickel1993efficient,tsiatis2007semiparametric}. They suggest sample splitting which together with the above Neyman orthogonality leads to ignorable remainder terms even when machine learning nuisance estimators are converging relatively slowly.
Other works have considered a Neyman orthogonal estimating equation of a $l^2$-continuous functional of a conditional expectation \citep{chernozhukovRiesz, bradic2019sparsity}. In this setting, both the conditional expectation and the Reisz representer of the functional, which is the inverse propensity score in the case of $\tau$, must be estimated. However, in the latter case, the estimation of the inverse propensity score is performed differently compared to \citet{chernozhukov2018double}; i.e., using the equation which characterize the nuisance model as a Riesz representer.

The notion of double robustness of an estimator has been widely used to indicate that an estimator is consistent if at least one of the nuisance models is estimated consistently, not necessarily both \citep{robins1994estimation, bang2005doubly}. This property has also been called parametric double robustness \citep{kennedy2016semiparametric}. As mentioned in the previous section, 
\citet[][Section 4]{van2014targeted}  was the first work which addressed what they called double robust statistical inference (also called model double robustness in \citeauthor{smucler2019unifying}, \citeyear{smucler2019unifying}), which indicates that inference on the parameter of interest requires consistent estimation of only one of the nuisance models. The AIPW estimator mentioned in the previous section based on regularized maximum likelihood nuisance estimators has the nonparametric double robust property \citep{kennedy2016semiparametric}, also called rate double robustness in \citet{smucler2019unifying}, whereby weak consistency of the nuisance models is required for uniformly valid inference, but slower convergence in estimating the propensity score can be bought out by faster convergence in the outcome models, and vice versa (\ref{productrate.eq}). However, this AIPW estimator does not yield double robust statistical inference. Alternative loss functions have been considered in the estimation of nuisance models, which endow the AIPW estimator double robust statistical inference
\citep[e.g.][]{avagyan2017honest, tan2020model, ning2020robust, bradic2019sparsity}.
\citet{smucler2019unifying} extends this area of work by  generalizing the property of double robust inference to the estimation of all parameters that belong to what is called the class of bilinear influence function (BIF) functionals. These estimators are specific to sparse settings and employ $l_1$-regularized estimators. Moreover, they consider a sparsity condition even for the limit of a possibly inconsistent nuisance model estimator.

The BIF class includes important causal parameters such as the average treatment effect and the average treatment effect among treated and cover the classes of parameters studied in \citet{chernozhukovRiesz} and \citet{robins2008higher}. However, the parameters who enjoy a uniform valid inference are not limited to this class. For example, the continuously differentiable functions of the functionals with bilinear influence function do not belong to the class while validity of inference for those can be directly shown by the delta method \citep{smucler2019unifying}. A different technique has been used for constructing a valid confidence interval for a parameter outside the BIF class, a conditional average treatment effect, which is to invert a chi squared distributed double robust test statistic \citep{dukes2020inference}. 

Most of the above literature on uniformly valid causal inference is concerned with high-dimensional settings, which can arise both because of a large set of covariates is available, but also because functions/transformations of these covariates are considered in generalized linear (in the parameters) nuisance models with a sparsity property. However, alternative regularity conditions, e.g. smoothness, may be considered attractive. For example, neural networks can be used for smooth nuisance functions belonging to Sobolev spaces \citep{farrell2018deep}. Some of the results in the above-mentioned papers are not specific to $l_1$-regularized estimators and apply to any well-behaved nonparametric nuisance model estimation \citep{van2006targeted, farrell2015robust, chernozhukov2018double}. In \citet{van2006targeted} the choice of a single nonparametric estimator is not considered to be done apriori, but a data-adaptive cross-validated combination of a set of estimators (super learner, ensemble learner) was suggested. Finally, more recently, \citet{cui2019selective} suggest two novel selection criteria, where the main focus is on getting smaller bias in the estimation of the target parameter instead of the nuisance ones. 

While the literature has focused on situations where the parameter (a causal effect) is of low dimension, \citet{semenova2020debiased} recently addresses situations, where the Neyman orthogonal property can be applied to obtain uniformly valid results in nonparametric situations, i.e. where the parameter of interest is of infinite dimension. Examples include causal effects conditional on a continuous covariate and causal effects of a continuous valued treatment.

\section{Cost of uniformity and a double-selection outcome regression estimator}\label{ourestimator}

One essential component of estimators \textcolor{blue}{with uniformly valid inference reviewed} in the latter section is that if there are instruments $-$ here covariates that explain $T$ although they are not related to $Y$ conditional on the other covariates included in $\bm{X}$ $-$ these may be part of the selected set of covariates. This is also obviously true for propensity score centered methods; see, e.g., \citet{shortreed2017outcome}.
This is unfortunate because the semiparametric efficiency bound for the average treatment effect is lower if we have knowledge on which variables are instruments \citep{hahn1998role,hahn2004functional,de2011covariate,RotnitzkySmucler:2020,tang2020outcome}. The variance inflation due to instruments can be severe and this has been reported in the literature numerous times, see, \citet{Schnitzer2016} and references therein. 

In the sequel, we provide a discussion and results which shed new light on this issue by presenting an estimation strategy which seems to yield a compromise between the estimators for which we have a uniformly valid asymptotic distribution (using instruments) and superefficient estimators, where irrelevant instruments are selected away by using the data.

For simplicity, consider as parameter of interest  $\tau_1 = E(Y(1))$. However, the results for $\tau_0 = E(Y(0))$ and thereby $\tau=\tau_1-\tau_0$ are analogous. Let $\bm{x} = [\bm{x}_1 , \cdots, \bm{x}_n]'$, $\bm{y}= [y_1 , \cdots , y_n]'$,  superscript $T$ denotes subsetting rows that correspond to treated individuals and subscript $S$ denotes subsetting columns using the set $S$. Then,  $P^T_S=\bm{x}^T_{S} (\bm{x}^{T \prime}_{S} \bm{x}^{T }_{S})^{-1} \bm{x}^{T \prime}_{S}$ is the projection matrix onto the space spanned by $\bm{x}^T_{S}$, while $\tilde{P}^T_S=\bm{x}_{S} (\bm{x}^{T \prime}_{S} \bm{x}^{T }_{S})^{-1} \bm{x}^{T \prime}_{S}$ is the matrix used in predicting $Y(1)$ for all $n$ individuals. Given a selected covariate set $S$, we define the post-selection outcome regression (OR) estimator as
 \begin{equation*} \label{or}
    \hat{\tau}_{1,OR}(S) = E_n[(\tilde{P}^T_S \bm{y}^T)_i]
    = \dfrac{1}{n} \Sigma_{i=1}^n [\bm{x}_{S,i} (\bm{x}^{T \prime}_{S} \bm{x}^{T }_{S})^{-1} \bm{x}^{T \prime}_{S} \bm{y}^T  ] 
    :=\dfrac{1}{n} \Sigma_{i=1}^n [ \hat{m}_1(\bm{x}_{{S},i}) ]
    .\end{equation*}
A classical post-selection OR estimator is $\hat{\tau}_{1,OR_{XY}}=\hat{\tau}_{1,OR}(X_Y)$, where $X_Y$ is a set of covariates derived from the $y-\bm{x}$ association using any covariate selection strategy; for instance, the set of covariates that corresponds to nonzero coefficients in a fitted lasso regression of $y$ versus $\bm{x}$ may be considered.
Instead, we study here theoretically the post-double-selection OR estimator
\begin{equation}
    \label{postdoubleor}\hat{\tau}_{1,OR_{DS}}=\hat{\tau}_{1,OR}(X_Y \cup X_T),
\end{equation}
 where $X_T$ is derived by fitting the  $t-\bm{x}$ association using lasso or any other covariate selection strategy. Note that this can be considered as a generalization of \cite{belloni2014inference} early estimator for the homogeneous case presented in Section \ref{homog.sec} to the general situation of an heterogeneous treatment effect. Although this estimator is mentioned in the simulation experiments run in \cite{athey2018approximate}, no theoretical results are available in the literature up to our knowledge.

Estimator (\ref{postdoubleor}) is not asymptotically linear and cannot be shown to yield uniformly valid inference as was the case for the post-selection double robust estimators described in Section \ref{hetero.sec}.
 However, uniform fast rate of decay of the bias can be guaranteed. The conditions used below are of the type used to show uniform validity. In particular, a product convergence rate condition is used. 
\begin{theorem}\label{theorem1} Suppose
\begin{enumerate}[label={\upshape(\roman*)}]
\item$E_n[(1 - t_i a(\bm{x}_i))( \tilde{P}^T_S m_1(\bm{x}^T) - {m}_1(\bm{x}) )_{i} ]= o_{P_n}(n^{-v_1}),$
\vspace{4pt}
\item
    $ E_{n_t}[((\mathbb{1}_{n_t\times n_t} - P^T_S) a(\bm{x}^T) )_i^2]^{1/2} E_{n_t}[((\mathbb{1}_{n_t\times n_t} - P^T_S) m_1(\bm{x}^T) )_i^2]^{1/2} = o_p(n^{-v_2}),$
\end{enumerate} 
 where $\mathbb{1}_{n_t\times n_t}$ is the identity matrix of size $n_t$, $a(\bm{X}) = 1/\mathbb{P}(T=1|\bm{X})$. Let $v = \min (v_1,v_2)$. Then, 
\begin{flalign*}
  \text{Bias} ( \hat{\tau}_{1,OR}) = E(E_n[\hat{m}_1(\bm{x}_{S,i})] -\tau_{1}) = o(n^{-v}).
\end{flalign*}

\end{theorem}
\noindent The first condition (i)  requires that the order of the scaled error for all the individuals is equal to the order of the scaled error on the treated individuals weighted by inverse propensity scores, similar to \citet[][Assumption 3(c)]{farrell2015robust}. The second condition (ii) is a multiplicative rate condition similar to the multiplicative rate condition in \citet[][Assumption 3(b)]{farrell2015robust}. However, notice that here the rate must be fulfilled using the same set of covariates in both nuisance models. This necessitates doing double selection to get the double robustness property in terms of bias. The proof of Theorem \ref{theorem1} can be found in Appendix \ref{A1}.

To illustrate how a double-selection procedure can benefit in terms of bias, suppose
\begin{equation}\label{ourmodely}
m_1(\bm{X},\eta_1)= \eta'_{1}\bm{X}+ R_Y,
\end{equation}
and
\begin{equation}\label{ourmodelt}
a(\bm{X},\gamma)= \gamma'\bm{X} + R_T.
\end{equation}
Here, we consider a linear model for the inverse of propensity score model \cite[e.g.][]{imbens2005mean}, and $R_Y$ and $R_t$ are approximation errors of the true models with sparse coefficients $\eta_1$ and $\gamma$ in the outcome and inverse propensity score models, respectively. The following corollary is a direct consequence of Theorem \ref{theorem1} and the asymptotic results for lasso regression in \citet[][Section 6]{farrell2015robust}.
\begin{corollary}\label{corollary1}
Suppose that models (\ref{ourmodely}) and (\ref{ourmodelt}), hold, where the conditions (\ref{errorY}), (\ref{errorT}), (\ref{consistencysparsity}) and (\ref{ratesparsity}) are fulfilled. Moreover, assume the regularity conditions in \citet[][Corollary 5 and Appendix F.3]{farrell2015robust}. Consider $\hat{\tau}_{1,OR_{DS}}$ in (\ref{postdoubleor}) where $X_Y$ and $X_T$ are estimated using lasso regression of the observed potential outcome (for treated individuals) and the inverse propensity score on $X$, respectively. 
Then,
$$
 \text{Bias} (\hat{\tau}_{1,OR_{DS}}) = o(n^{-1/2}).
$$
\end{corollary}
The above result shows that root-$n$ decay of the bias can be derived uniformly over a set of DGPs. In this sense, the double-selection OR estimator may be seen as a compromise between single selection estimators (superefficient and no uniformly decaying bias) and the double-selection estimators of Section \ref{uniformlyvalidinference} (with available uniformly valid asymptotic distribution). The implications in terms of finite sample behaviour are studied below with Monte Carlo experiments.

In practice, the inverse propensity score $a$ is not observed and (\ref{ourmodelt}) cannot be fitted directly. Instead, we suggest to fit a lasso logistic regression for the propensity score to retrieve the relevant covariates.

\section{Simulation study}\label{simulation}
The aim of this simulation study is to illustrate the above theoretical discussion on the cost and benefits of uniformly valid inference. While many simulations studies are available in the literature reviewed above, their focus is either to show the necessity of using uniformly valid procedures in order to avoid regularization bias, or to illustrate the variance inflation due to the use of instruments in the propensity score compared to superefficient procedures. Here, we contrast these two aspects by considering both uniformly valid and superefficient post-selection strategies, as well as the double-selection outcome regression estimator (\ref{postdoubleor}). 

\subsection{Simulation design}

We use 500 replicates in all situations, and consider sample sizes $n$ = 500, 1000, 1500, 2000. The covariate vector $\bm{X}$ is generated from a multivariate normal distribution with zero mean and identity covariance matrix. The dimension $p$ of the covariate vector equals $n$. Results for low dimensional settings, $p<<n$, portray a similar general picture and are available from the authors upon request. Data generation and all computations are performed with the software R \citep{R}.

\subsubsection{High-dimensional setting}
We consider the following models

$$
Y(0) = m_0(\bm{X}) + \epsilon_0 = 1+ \eta_0'\bm{X} +\epsilon_0,
$$
$$
Y(1) = m_1(\bm{X}) + \epsilon_1= 2+  \eta_1'\bm{X} +\epsilon_1,
$$
$$
\mathbb{P}(T=1|\bm{X}) = e(\bm{X}) = \text{expit}(\gamma'\bm{X}),
$$

indexed by the parameter vectors
$$
\eta_0 =  \frac{k}{2} \cdot (1,1/2,1/3,1/4,1/5,0,0,0,0,0,
      1,1/2,1/3,1/4,1/5,0, \cdots,0),
$$
$$
\eta_1 = k \cdot (1,1/2,1/3,1/4,1/5,0,0,0,0,0,
      1,1/2,1/3,1/4,1/5,0, \cdots,0),
$$
$$
\gamma = (1,1/2,1/3,1/4,1/5,1,1,1,1,1,0, \cdots,0),
 $$ 
where $k\in\{0.1, 0.4, 0.8, 1.2\}$ and the error terms, $\epsilon_t,\,\, t= 0, 1$, are generated from a normal distribution with $E(\epsilon_t|\bm{X})=0$ and $Var(\epsilon_t|\bm{X}) = (1+p)^{-1}(1+\iota'\cdot \bm{X}^2)$ where $\iota$ is the vector of ones. The parameter $k$ determines the strength of the association between the outcome and the covariates.
 
\subsubsection{Post-selection estimators of $\tau$}
We use lasso as implemented in the R package \texttt{hdm} \citep{chernozhukov2016hdm} to estimate the nuisance models $m_t(\bm{X})$ and $e(\bm{X})$. We denote the sets of variables which corresponds to nonzero coefficients in the estimated sparse linear outcome models and logistic propensity score by $X_Y = X_{Y_0}\cup X_{Y_1}$, the union of the sets estimated by each of the two potential outcome models, and $X_T$, respectively. The lasso penalty parameter is selected as $\lambda = 2.2\sqrt{n}\Phi^{-1}([\mbox{log}(n) - 0.1][2p\mbox{log}(n)]^{-1})$ \citep{BCCH:12}.

With the purpose of estimating $\tau$, using different combinations of the above covariate sets, we compare two versions of the OR estimator, three versions of the doubly robust AIPW estimator and one version of the doubly robust TMLE estimator. 
Specifically, OR$_{X_Y}$ and OR$_{DS}$ use $X_Y$ and the union $X_T\cup X_Y$ in the outcome models refitting steps, respectively. AIPW$_{X_Y}$ uses $X_Y$ in the propensity score and outcome model refitting steps, AIPW$_{DS}$ uses the union $X_T \cup X_Y$ in both the propensity score and outcome models  \citep[as in][Section 5]{belloni2014inference} and AIPW$_{X_Y, X_T}$ uses $X_T$ in the propensity score model and $X_Y$ in the outcome models \citep[as in][]{farrell2015robust}. TMLE$_{X_Y, X_T}$ uses $X_T$ in the propensity score model and $X_Y$ in the outcome models. 

Given the same consistency and product rate conditions on the initial estimators of nuisance models as in  \citet{farrell2015robust}, TMLE$_{X_Y, X_T}$ has the same uniformly valid asymptotic distribution as AIPW$_{X_Y, X_T}$ \citep[][Chapter 27]{vdLR:11}. Hence, if the refitted models used in AIPW$_{X_Y, X_T}$ are used as initial model estimates in TMLE$_{X_Y, X_T}$ we would expect similar results for large samples. In summary, AIPW$_{X_Y, X_T}$, AIPW$_{DS}$ and TMLE$_{X_Y, X_T}$ have uniformly valid asymptotic distributions, OR$_{X_Y}$ and  AIPW$_{X_Y}$ have no such uniform validity, select away instruments and are superefficient, while OR$_{DS}$ has uniformly decaying bias (Corollary \ref{corollary1}). For TMLE, we use the R package \texttt{tmle} \citep{tmlepackage} and do not truncate the estimated propensity scores, i.e., \texttt{gbound = c(0,1)}. For the other estimators, we use own written R code as well as the R package \texttt{ui} \citep{genbackdeluna19} for the variances.

\subsection{Results}
From 500 replicates we compute empirical biases, standard errors, root mean squared errors (RMSE), and empirical coverages. We also compute mean estimated standard errors. Table \ref{tab:high} presents results for $k=0.4$. Results for all values of $k$ are given for bias, RMSE and coverages in the Appendix, Tables \ref{tab.bias}-\ref{tab.ec}.

We see that estimators selecting away instruments (OR$_{X_Y}$ and AIPW$_{X_Y}$) have lower Monte Carlo standard error but at the cost of larger bias and poor empirical coverages (clearly lower than nominal level). As expected the other estimators, which all have uniformly decaying bias, show low bias, but at the cost of larger standard error. This cost is, however, smallest for OR$_{DS}$, and AIPW$_{X_Y, X_T}$ , AIPW$_{DS}$ and TMLE$_{X_Y, X_T}$  have standard errors roughly up to five times as large as OR$_{X_Y}$, while for OR$_{DS}$ the increase in standard error is not as severe (1.22 times the standard error of OR$_{X_Y}$). All the low-bias estimators have good empirical coverages, at least for sample sizes 1000 and higher, although we do not have such theoretical guarantee for OR$_{DS}$.

On a side note, we observe that the estimated standard errors of AIPW$_{X_Y, X_T}$, AIPW$_{DS}$ and TMLE$_{X_Y, X_T}$ are distinctively smaller than the Monte Carlo standard errors, but the difference is reduced when increasing sample size.  For OR$_{DS}$ no such underestimation of the variance is observed.

\begin{table}[!htb]
\vspace{4mm}
\centering
\caption{Results of 500 simulation replicates for estimators of $\tau$, for varying sample sizes $n$, number of covariates $p=n$ and $k=0.4$. RMSE, root mean-squared error; Bias; SE, Monte Carlo standard deviation; ESE, estimated standard error (influence curve based estimates, ignoring the variability in the selection step); CP, empirical coverage probability of 95\% confidence intervals.}
\vskip 1mm \vskip 0pt
\label{tab:high}
\begin{tabular}{llrrrrr}
\hline
  \hline
\multicolumn{1}{c}{\rule{0pt}{3ex}$n$}& \multicolumn{1}{l}{Estimator} &\multicolumn{1}{c}{RMSE} & \multicolumn{1}{c}{Bias} & \multicolumn{1}{c}{SE} &\multicolumn{1}{c}{ESE} & \multicolumn{1}{c}{CP}  \\ 
\hline
\rule{0pt}{3ex}500 
&\1 & 0.26 & -0.19 & 0.18 & 0.16 & 0.74 \\ 
&\2 & 0.26 & -0.19 & 0.18 & 0.16 & 0.74 \\ 
&\3 & 0.22 & -0.05 & 0.22 & 0.21 & 0.94 \\ 
&\4 & 0.95 & -0.06 & 0.95 & 0.35 & 0.93 \\ 
&\5 & 0.88 & -0.06 & 0.88 & 0.34 & 0.91 \\ 
&\6 & 0.29 & -0.08 & 0.28 & 0.24 & 0.91 \\ 

\rule{0pt}{3ex}1000 
&\1 & 0.19 & -0.14 & 0.14 & 0.11 & 0.72 \\ 
&\2 & 0.19 & -0.14 & 0.14 & 0.11 & 0.73 \\ 
&\3 & 0.16 & -0.02 & 0.16 & 0.15 & 0.94 \\ 
&\4 & 0.38 & -0.00 & 0.38 & 0.25 & 0.97 \\ 
&\5 & 0.37 & -0.00 & 0.37 & 0.25 & 0.95 \\ 
&\6 & 0.21 & -0.04 & 0.20 & 0.18 & 0.95 \\ 

\rule{0pt}{3ex}1500 
&\1 & 0.15 & -0.11 & 0.11 & 0.09 & 0.76 \\ 
&\2 & 0.15 & -0.11 & 0.11 & 0.09 & 0.76 \\ 
&\3 & 0.12 & -0.02 & 0.12 & 0.12 & 0.95 \\ 
&\4 & 0.22 & -0.03 & 0.22 & 0.18 & 0.95 \\ 
&\5 & 0.24 & -0.02 & 0.24 & 0.18 & 0.95 \\ 
&\6 & 0.15 & -0.04 & 0.15 & 0.14 & 0.92 \\ 
  
\rule{0pt}{3ex}2000
&\1 & 0.13 & -0.09 & 0.10 & 0.08 & 0.80 \\ 
&\2 & 0.13 & -0.09 & 0.10 & 0.08 & 0.80 \\ 
&\3 & 0.10 & -0.01 & 0.10 & 0.10 & 0.95 \\ 
&\4 & 0.19 & -0.02 & 0.19 & 0.16 & 0.94 \\ 
&\5 & 0.19 & -0.02 & 0.19 & 0.16 & 0.93 \\ 
&\6 & 0.15 & -0.03 & 0.15 & 0.13 & 0.92 \\
   \hline
\end{tabular}
\end{table}

\section{Discussion}
To obtain an unbiased estimator of the average causal effect of a treatment, we need to control for all confounders. Knowing that some covariates are related to the outcome, but not to the treatment does not change the semiparametric efficiency bound. However, knowledge on the existence of instruments has implications on the asymptotic variance that can be achieved by unbiased estimators \citep{hahn2004functional,de2011covariate,RotnitzkySmucler:2020}. In practice, we typically do not have such apriori (to data) knowledge. Naively selecting away covariates (here instruments) using the data at hand, by, e.g., regularization yields inference that is not uniformly valid (may translate into large bias and incorrect coverage rates). On the other hand, using methods which yield uniformly valid inference may yield large variability compared to these naive (often superefficient) methods. In this paper, we have reviewed the literature on uniformly valid causal inference, and have discussed the costs and benefits of uniformly valid inference. The latter discussion has been illustrated by studying 
a double-selection outcome regression estimator, which is shown to be uniformly asymptotically unbiased under a product rate condition. This seems to translate into finite sample properties which are a compromise between uniformly valid and superefficient estimators. The good properties of the double-selection OR estimator may arguably be due to the designs considered in our simulations, for which the outcome regression models can be consistently fitted. Consistency of the outcome regression model is indeed assumed to show uniformly decaying bias of the double-selection estimator. This assumption is also used by estimators which have the nonparametric double robust property. On the other hand, procedures that yield double robust statistical inference \citep{van2014targeted,tan2020model,avagyan2017honest}, allow for one of the nuisance models to be inconsistent (converges to an incorrect function of the covariates). This alternative limit need, however, to be assumed sparse in the covariates, a rather high-level assumption \citep{smucler2019unifying}.

\section*{Acknowledgements}
We acknowledge comments from Tetiana Gorbach, Mohammad Ghasempour and two anonymous referees that have improved the paper. This work was supported by a grant from the Marianne and Marcus Wallenberg Foundation.

\clearpage
\appendix
\label{appendix}
\section{Appendix}
\subsection{Proof of Theorem 1}\label{A1}

The bias of post selection outcome regression estimator can be expressed as follows
\begin{flalign*}
    \text{Bias}(\hat{\tau}_{1,OR}) &= E(E_n[\hat{m}_1(\bm{x}_{S,i})] -\tau_{1})&  \\&= E  (E_n[\hat{m}_1(\bm{x}_{S,i}) -  {m}_1(\bm{x}_{i})] )  +  E(E_n [{m}_1(\bm{x}_{i}) ] - \tau_1)&\\ &= E  (E_n[\hat{m}_1(\bm{x}_{S,i}) -  {m}_1(\bm{x}_{i})] ),  &
\end{flalign*}
where the last equality follows by Assumption \ref{NUC}. To show that the scaled bias term is asymptotically negligible we show that it is negligible conditional on $S$ and $\{t_i,x_i\}_{i=1}^n$.
\begin{flalign*}
 n^v   \text{Bias}(\hat{\tau}_{1,OR}|S,\{t_i,x_i\}_{i=1}^n) &:= n^v E  (E_n[\hat{m}_1(\bm{x}_{S,i}) -  {m}_1(\bm{x}_{i})]|S,\{t_i,x_i\}_{i=1}^n ) &\\&=
n^v E_n[ \bm{x}_{S,i} (\bm{x}^{T \prime}_{S} \bm{x}^{T }_{S})^{-1} \bm{x}^{T \prime}_{S} m_1(\bm{x}^T) - {m}_1(\bm{x}_{i}) ] &\\&
\approx n^v E_n[ t_i a(\bm{x}_i)( \bm{x}_{S,i} (\bm{x}^{T \prime}_{S} \bm{x}^{T }_{S})^{-1} \bm{x}^{T \prime}_{S} m_1(\bm{x}^T) - {m}_1(\bm{x}_{i}) ) ]& \\
& =\dfrac{n^v}{n}  a(\bm{x}^T)'( P^T_S - \mathbb{1}_{n_t\times n_t}) {m}_1(\bm{x}^T)  &\\& =\dfrac{n^v}{n}  (( P^T_S - \mathbb{1}_{n_t\times n_t})a(\bm{x}^T))'( P^T_S - \mathbb{1}_{n_t\times n_t}) {m}_1(\bm{x}^T),  &
\end{flalign*}
 where
 $\approx$ means that both side have the same limits which holds by \ref{theorem1}{(i)}.
\begin{flalign*}| n^v   \text{Bias}(\hat{\tau}_{1,OR}|S,\{t_i,x_i\}_{i=1}^n)| &\leq
n^v E_{n_t}[((P^T_S - \mathbb{1}_{n_t\times n_t})a(\bm{x}^T))_i^2]^{1/2} &
\\&\ \ \ \ \times E_{n_t}[((P^T_S - \mathbb{1}_{n_t\times n_t})m_1(\bm{x}^T))_i^2]^{1/2}& \\&= o_{P_n}(1),&\end{flalign*}
where the last equality holds by \ref{theorem1}{(ii)}. The statement in the theorem follows by the above result on the order of decay of the conditional expectation and uniform integrability of this conditional expectation.

\subsection{Proof of Corollary 1}
By construction $X_Y \subseteq S$. Therefore by \citet[][Appendix F.3]{farrell2015robust}, (\ref{errorY}) and the sparsity condition (\ref{consistencysparsity})  we have $v_1 = 1/2$.
Moreover, using \citet[][Corollary 5]{farrell2015robust} , $X_Y \subseteq S$, $X_T \subseteq S$,  and conditions (\ref{errorY}), (\ref{errorT}) and (\ref{ratesparsity}) we have $v_2 = 1/2$.

\vskip 5.5cm

\subsection{Simulation results}

  \begin{table}[h!]
    \centering
    \caption{$\sqrt{n}$-bias based on 500 simulation replicates, for estimators of $\tau$, with varying sample sizes $n$ and values for $k$, and number of covariates $p=n$.}
       \label{tab.bias}
\begin{tabular}{llrrrr}
\hline
\hline
&&\multicolumn{2}{c}{$k$}\\
\cline{3-6}
\multicolumn{1}{c}{\rule{0pt}{3ex}$n$} & Estimator & 0.1&0.4&0.8&1.2\\
\hline
\rule{0pt}{3ex}
500
&\1 & -1.25 & -4.16 & -3.44 & -2.97 \\ 
&\2 & -1.25 & -4.17 & -3.44 & -2.96 \\ 
&\3 & -0.26 & -1.08 & -2.09 & -2.56 \\ 
&\4 & -0.42 & -1.26 & -2.57 & -3.12 \\ 
&\5 & -0.35 & -1.45 & -2.49 & -2.44 \\ 
&\6 & -0.45 & -1.79 & -2.55 & -2.84 \\
\rule{0pt}{3ex}
1000
&\1 & -1.81 & -4.32 & -3.14 & -2.28 \\ 
&\2 & -1.81 & -4.31 & -3.15 & -2.27 \\ 
&\3 & -0.17 & -0.76 & -1.40 & -1.63 \\ 
&\4 & 0.59 & -0.06 & -0.86 & -1.12 \\ 
&\5 & 0.56 & -0.03 & -0.88 & -1.08 \\ 
&\6 & -0.30 & -1.23 & -1.76 & -1.83 \\ 
\rule{0pt}{3ex}
1500
&\1 & -2.50 & -4.10 & -3.20 & -2.08 \\ 
&\2 & -2.50 & -4.09 & -3.21 & -2.08 \\ 
&\3 & -0.50 & -0.91 & -1.44 & -1.47 \\ 
&\4 & -0.66 & -1.02 & -1.36 & -1.55 \\ 
&\5 & -0.42 & -0.89 & -1.39 & -1.48 \\ 
&\6 & -0.75 & -1.38 & -1.71 & -1.59 \\ 
\rule{0pt}{3ex}
2000
&\1 & -2.83 & -3.91 & -2.91 & -1.56 \\ 
&\2 & -2.82 & -3.92 & -2.89 & -1.53 \\ 
&\3 & -0.29 & -0.60 & -0.91 & -0.87 \\ 
&\4 & -0.54 & -0.91 & -1.21 & -1.02 \\ 
&\5 & -0.44 & -0.81 & -1.08 & -0.85 \\ 
&\6 & -0.68 & -1.15 & -1.33 & -1.04 \\ 
 \hline
\end{tabular}
    \end{table}

  \begin{table}[t]
    \centering
    \caption{RMSE based on 500 simulation replicates, for estimators of $\tau$, with varying sample sizes $n$ and values for $k$, and number of covariates $p=n$.}
       \label{tab.rmse}
\begin{tabular}{llrrrr}
\hline
\hline
&&\multicolumn{2}{c}{$k$}\\
\cline{3-6}
\multicolumn{1}{c}{\rule{0pt}{3ex}$n$} & Estimator & 0.1&0.4&0.8&1.2\\
\hline
\rule{0pt}{3ex}
500
&\1 & 0.17 & 0.26 & 0.25 & 0.24 \\ 
&\2 & 0.17 & 0.26 & 0.25 & 0.24 \\ 
&\3 & 0.21 & 0.22 & 0.24 & 0.26 \\ 
&\4 & 0.85 & 0.95 & 0.95 & 0.88 \\ 
&\5 & 0.79 & 0.88 & 0.88 & 0.76 \\ 
&\6 & 0.28 & 0.29 & 0.31 & 0.32 \\ 
  \rule{0pt}{3ex}
1000
&\1 & 0.12 & 0.19 & 0.17 & 0.15 \\ 
&\2 & 0.12 & 0.19 & 0.17 & 0.15 \\ 
&\3 & 0.15 & 0.16 & 0.17 & 0.17 \\ 
&\4 & 0.37 & 0.38 & 0.37 & 0.38 \\ 
&\5 & 0.36 & 0.37 & 0.38 & 0.38 \\ 
&\6 & 0.20 & 0.21 & 0.22 & 0.22 \\
\rule{0pt}{3ex}
1500
&\1 & 0.11 & 0.15 & 0.13 & 0.11 \\ 
&\2 & 0.11 & 0.15 & 0.13 & 0.11 \\ 
&\3 & 0.12 & 0.12 & 0.13 & 0.13 \\ 
&\4 & 0.22 & 0.22 & 0.24 & 0.24 \\ 
&\5 & 0.22 & 0.24 & 0.25 & 0.25 \\ 
&\6 & 0.15 & 0.15 & 0.16 & 0.16 \\ 
\rule{0pt}{3ex}
2000
&\1 & 0.10 & 0.13 & 0.11 & 0.10 \\ 
&\2 & 0.10 & 0.13 & 0.11 & 0.10 \\ 
&\3 & 0.10 & 0.10 & 0.11 & 0.11 \\ 
&\4 & 0.18 & 0.19 & 0.19 & 0.19 \\ 
&\5 & 0.18 & 0.19 & 0.19 & 0.19 \\ 
&\6 & 0.15 & 0.15 & 0.15 & 0.15 \\
 \hline
\end{tabular}
    \end{table}

 \begin{table}[t]
    \centering
    \caption{Empirical coverage probability of 95\% confidence intervals based on 500 simulation replicates, for estimators of $\tau$, with varying sample sizes $n$ and values for $k$, and number of covariates $p=n$..
}
       \label{tab.ec}
\begin{tabular}{llrrrr}
\hline
\hline
&&\multicolumn{2}{c}{$k$}\\
\cline{3-6}
\multicolumn{1}{c}{\rule{0pt}{3ex}$n$} & Estimator & 0.1&0.4&0.8&1.2\\
\hline
\rule{0pt}{3ex}
500
&\1 & 0.94 & 0.74 & 0.83 & 0.86 \\ 
&\2 & 0.94 & 0.74 & 0.82 & 0.87 \\ 
&\3 & 0.93 & 0.94 & 0.92 & 0.92 \\ 
&\4 & 0.95 & 0.93 & 0.92 & 0.91 \\ 
&\5 & 0.92 & 0.91 & 0.90 & 0.90 \\ 
&\6 & 0.93 & 0.91 & 0.89 & 0.90 \\ 
\rule{0pt}{3ex}
1000
&\1 & 0.93 & 0.72 & 0.84 & 0.89 \\ 
&\2 & 0.93 & 0.73 & 0.84 & 0.90 \\ 
&\3 & 0.93 & 0.94 & 0.92 & 0.92 \\ 
&\4 & 0.97 & 0.97 & 0.95 & 0.95 \\ 
&\5 & 0.96 & 0.95 & 0.94 & 0.93 \\ 
&\6 & 0.95 & 0.95 & 0.94 & 0.94 \\ 
\rule{0pt}{3ex}
1500
&\1 & 0.90 & 0.76 & 0.83 & 0.91 \\ 
&\2 & 0.90 & 0.76 & 0.84 & 0.92 \\ 
&\3 & 0.96 & 0.95 & 0.96 & 0.95 \\ 
&\4 & 0.96 & 0.95 & 0.95 & 0.94 \\ 
&\5 & 0.95 & 0.95 & 0.94 & 0.93 \\ 
&\6 & 0.95 & 0.92 & 0.93 & 0.92 \\
\rule{0pt}{3ex}
2000
&\1 & 0.88 & 0.80 & 0.86 & 0.92 \\ 
&\2 & 0.88 & 0.80 & 0.87 & 0.92 \\ 
&\3 & 0.95 & 0.95 & 0.94 & 0.95 \\ 
&\4 & 0.94 & 0.94 & 0.94 & 0.95 \\ 
&\5 & 0.94 & 0.93 & 0.94 & 0.94 \\ 
&\6 & 0.92 & 0.92 & 0.92 & 0.93 \\ 
 \hline
\end{tabular}
    \end{table}

\clearpage

\bibliography{biblioDR}

\end{document}